\documentclass[12pt,leqno]{article}
\usepackage{amsmath,amssymb,amsfonts,amsthm,latexsym,amstext,amscd}
\usepackage[english]{babel}
\usepackage{fancyhdr}

\newtheorem{theorem}{Theorem}
\newtheorem{corollary}{Corollary}
\newtheorem{lemma}{Lemma}

\newtheorem*{theo}{Theorem}

\theoremstyle{definition}

\newtheorem{remark}{Remark}

\newtheorem*{rema}{Remark}

\newcommand{\beq}{\begin{equation}}
\newcommand{\eeq}{\end{equation}}

\def\ve{\varepsilon }
\begin{document}

\title
{On the distribution of gaps between consecutive primes}

\author
{J\'anos Pintz}

\date{}

\numberwithin{equation}{section}


\maketitle

\section{Introduction}
\label{sec:1}

The recent dramatic new developments in the study of bounded gaps between primes, reached by Zhang \cite{Zha}, Maynard \cite{May1} and Tao \cite{Pol8B} made other old conjectures about the distribution of primegaps
\beq
\label{eq:1.1}
d_n = p_{n + 1} - p_n, \quad \mathcal P = \{p_i\}_{i = 1}^\infty \ \text{ the set of all primes}
\eeq
accessible.
One of the most interesting such conjectures was formulated in 1954 by Erd\H{o}s \cite{Erd2} as follows.
Let $J$ denote the set of limit points of $d_n / \log n$, i.e.\
\beq
\label{eq:1.2}
J = \left\{\frac{d_n}{\log n}\right\}'.
\eeq
Then $J = [0, \infty]$.

While Westzynthius \cite{Wes} proved already in 1931 the relation
\beq
\label{eq:1.3}
\limsup_{n \to \infty} \frac{d_n}{\log n} = \infty \ \ \text{ i.e. }\ \ \infty \in J,
\eeq
no finite limit point was known until 2005, when in a joint work of Goldston, Y{\i}ld{\i}r{\i}m and the author \cite{GPY1} it was shown that
\beq
\label{eq:1.4}
\liminf_{n \to \infty} \frac{d_n}{\log n} = 0 \ \ \text{ i.e. }\ \ 0 \in J.
\eeq
On the other hand, Erd\H{o}s \cite{Erd2} and Ricci \cite{Ric} proved simultaneously and independently about 60 years ago that $J$ has positive Lebesgue measure.

In a recent work D. Banks, T. Freiberg and J. Maynard showed \cite{BFM} that more than 2\% of all nonnegative real numbers belong to~$J$.

The author has shown \cite{Pin3} that for any $f(n) \leqslant \log n$, $f(n) \nearrow \infty$ (i.e.\ $f(n) \to \infty$, $f(n)$ is monotonically increasing) satisfying for any $\varepsilon$
\beq
\label{eq:1.5}
(1 - \ve) f(N) \leqslant f(n) \leqslant (1 + \ve) f(N) \ \text{ if } \ n \in [N, 2N], \ N > N_0(\ve)
\eeq
we have an ineffective constant $c_f$ such that
\beq
\label{eq:1.6}
[0, c_f] \subset J_f := \left\{ \frac{d_n}{f(n)}\right\}'.
\eeq

Although $\log n$ is the average value of $d_n$, improving the result \eqref{eq:1.3} of Westzynthius, Erd\H{o}s \cite{Erd1} in 1935 and three years later Rankin \cite{Ran1} proved stronger results about large gaps between consecutive primes.
The 76-year-old result of Rankin, the estimate ($\log_\nu n$ denotes the $\nu$-fold iterated logarithmic function)
\beq
\label{eq:1.7}
\limsup_{n \to \infty} \frac{d_n/\log n}{g(n)} \geqslant C_0, \ \ \ g(x) = \frac{\log_2 x \log_4 x}{(\log_3 x)^2}
\eeq
was apart from the value of the constant $C_0$ still until August 2014 the best known lower estimate for large values of $d_n$.
(The original value $C_0 = 1/3$ of Rankin was improved in four steps, finally to $C_0 = 2e^\gamma$ by the author \cite{Pin2}.)
Then in two days two different new proofs appeared by Ford--Green--Konjagin--Tao \cite{FGKT} and Maynard \cite{May2} in the arXiv, proving Erd\H{o}s's famous USD~10,000 conjecture according to which \eqref{eq:1.7} holds with an arbitrarily large constant~$C_0$.

This raises the question whether the relation \eqref{eq:1.6} can be improved to functions of type $f(x) = \omega(x) \log x$ with $\omega(x) \to \infty$ and whether perhaps even $\omega(x) = c_1 g(x)$ can be reached with some absolute constant~$c_1$, or, following the mentioned new developments, with an arbitrarily large $c_1$ as well.

Another question is whether for some function $f(n)$ we can reach
\beq
\label{eq:1.8}
[0, \infty] = J_f,
\eeq
i.e.\ the original conjecture of Erd\H{o}s with the function $f(n)$ in place of $\log n$.

Using our notation \eqref{eq:1.1}, \eqref{eq:1.5}--\eqref{eq:1.7}, we will show the following results, which, although do not show the original conjecture $J_{\log n} = [0, \infty]$ of Erd\H{o}s, but in several aspects approximate it and in other aspects they go even further.

Since the first version of this paper was written before the groundbreaking works \cite{FGKT} and \cite{May2} we will
present the formulation and proofs of the original version of our results in the Introduction and Sections~\ref{sec:uj3}--\ref{sec:uj5}, while the formulation of the improved stronger versions appear in Section~\ref{sec:uj2} and the needed changes in the proofs in Section~\ref{sec:uj6}, in this case the changes refer to the mentioned work of Maynard \cite{May2}.

\begin{theorem}
\label{th:1}
There exists an absolute constant $c_0$ such that for any function $f(x) \nearrow \infty$, satisfying \eqref{eq:1.5} and
\beq
\label{eq:1.9}
f(x) \leqslant c_0 g(x) \log x
\eeq
we have with a suitable (ineffective) constant $c_f$
\beq
\label{eq:1.10}
[0, c_f] \subset J_f := \left\{\frac{d_n}{f(n)}\right\}'.
\eeq
\end{theorem}

\begin{theorem}
\label{th:2}
Let us consider a sequence of functions $\left\{f_i(x)\right\}_{i = 1}^\infty$ satisfying \eqref{eq:1.5}, \eqref{eq:1.9}, $f(x) \nearrow \infty$ and
\beq
\label{eq:1.11}
\frac{f_{i + 1}(x)}{f_i(x)} \to \infty \ \text{ as } \ x \to \infty \ \text{ for every } i.
\eeq
Then apart from at most $98$ functions $f_i(x)$ we have
\beq
\label{eq:1.12}
[0, \infty] = J_{f_i} = \left\{\frac{d_n}{f_i(n)}\right\}'.
\eeq
\end{theorem}

Answering a question raised in a recent work of Banks, Freiberg and Maynard \cite{BFM} we show that the method of \cite{BFM} works
also if we normalize the primegaps in place of $\log n$ with any function not exceeding the Erd\H{o}s--Rankin function.

\begin{theorem}
\label{th:3}
Suppose $f(x) \nearrow \infty$ and satisfies \eqref{eq:1.5} and \eqref{eq:1.9}.
Then for any sequence of $k \geqslant 50$ non-negative real numbers
$\beta_1 < \beta_2 < \dots < \beta_k$ at least one of the numbers $\{\beta_j - \beta_i;\, 1 \leqslant i < j \leqslant k\}$ belongs to $J_f$.
Consequently more than $2\%$ of all non-negative real numbers belong to~$J_f$.
\end{theorem}

As a by-product the method also gives a different new proof for the following result of Helmut Maier \cite{Mai1} proved in 1981 by his famous matrix method:

\begin{theorem}
\label{th:4}
For any natural number $m$ we have with the notation \eqref{eq:1.7}
\beq
\label{eq:1.13}
\limsup_{h\to \infty} \frac{\min (d_{n + 1}, \dots, d_{n + m})}{g(n) \log n} > 0.
\eeq
\end{theorem}

An immediate corollary of Theorem~\ref{th:2} is the following

\begin{corollary}
\label{cor:1}
Let $\eta(x) \to 0$ be an arbitrary function.
If $\eta(x) F(x) \nearrow \infty$, $F(x) \nearrow \infty$, both functions $F(x)$ and $\eta(x) F(x)$ satisfy \eqref{eq:1.5} and \eqref{eq:1.9}, then we have a function $f(x) \nearrow \infty$,
\beq
\label{eq:1.14}
\eta(x)  F(x) \leqslant f(x) \leqslant F(x),
\eeq
for which
\beq
\label{eq:1.15}
[0, \infty] = J_f := \left\{\frac{d_n}{f(n)}\right\}'.
\eeq
\end{corollary}

This means that although we can not show Erd\H{o}s's conjecture for the natural normalizing function $\log n$, changing it a little bit, it will be already true for some function $\xi(n) \log n$, where $\xi(n)$ tends to $0$ (or alternatively we can require $\xi(n)\to \infty$) arbitrarily slowly (even if this is not a natural normalization).

\section{Stronger forms of Theorems~\ref{th:1}--\ref{th:4} and Corollary~\ref{cor:1}}
\label{sec:uj2}

We will use in the formulation and proof of our results the work of J. Maynard \cite{May2} which implicitly defines an unspecified but actually explicitly calculable $\omega_0(x)$ function with the property
\beq
\label{eq:uj2.1}
\lim_{x \to \infty} \omega_0(x) = \infty,
\eeq
such that defining (cf.\ \eqref{eq:1.7})
\beq
\label{eq:uj2.2}
g_0(x) = \omega_0(x)g(x) = \omega_0(x) \frac{\log_2 x \log_4 x}{(\log_3 x)^2},
\eeq
the result \eqref{eq:1.7} holds with $g_0(n)$ in place of $g(n)$.

J. Maynard further mentions in the Remark at the end of his paper that he hopes to obtain his result with $\omega_0(x) = (\log_2 x)^{1 + o(1)}$ which would be the limit of the Erd\H{o}s--Rankin method.

We mention that such an improvement would almost surely lead to an improvement of our results too, since the mentioned idea (to show the same results with a uniformity in the variable $k$ of \cite{May2} for $k$ as large as $k \asymp (\log x)^\alpha$) would leave the structure of the proof unchanged.

We will denote by Theorems 1'--4' and Corollary 1' the stronger versions of Theorems~\ref{th:1}--\ref{th:4} and Corollary~\ref{cor:1}.
They are the following ($g_0(x)$ is defined in \eqref{eq:uj2.2}).

\smallskip
\noindent
{\bf Theorem 1'.}
{\it Theorem~\ref{th:1} holds with $g(x)$ replaced by $g_0(x)$ in \eqref{eq:1.9}.}

\smallskip
\noindent
{\bf Theorem 2'.}
{\it Theorem~\ref{th:2} holds with $g(x)$ replaced by $g_0(x)$ in \eqref{eq:1.9}.}

\smallskip
\noindent
{\bf Theorem 3'.}
{\it Theorem~\ref{th:3} holds with $g(x)$ replaced by $g_0(x)$ in \eqref{eq:1.9}.}

\smallskip
\noindent
{\bf Theorem 4'.}
{\it Theorem~\ref{th:4} holds with $g(x)$ replaced by $g_0(x)$ in \eqref{eq:1.9}.}

\smallskip
\noindent
{\bf Corollary 1'.}
{\it Corollary~\ref{th:1} holds with $g(x)$ replaced by $g_0(x)$ in \eqref{eq:1.9}.}


\section{The Maynard--Tao theorem}
\label{sec:uj3}

We call $\mathcal H_m = \{h_1, \dots, h_m\}$ an admissible $m$-tuple if $0 \leqslant h_1 < \dots < h_m$ and $\mathcal H_m$ does not occupy all residue classes $\text{\rm mod }p$ for any prime~$p$.

Further, we recall the Theorem of Landau--Page (see \cite[p.~95]{Dav}).

\begin{theo}
If $c_1$ is a suitable positive constant, $N$ arbitrary, there is at most one primitive character $\chi$ to a modulus $r \leqslant N$: for which $L(s, \chi)$ has a real zero $\beta$ satisfying
\beq
\label{eq:2.1}
\beta > 1 - \frac{c_1}{\log N}.
\eeq
\end{theo}

Such an exceptional character $\chi$ must be real, which means also that its conductor $r$ is squarefree, apart from the possibility that the prime $2$ appears in the factorization of $r$ with an exponent $2$ or $3$.
We have also
\beq
\label{eq:2.2}
\beta < 1 - \frac{c_2}{\sqrt{r}\log^2 r} \ \text{ \cite[p.~96]{Dav}}, \ \
\beta < 1 - \frac{c_3}{\sqrt{r}} \ \text{ \cite{Pin1, GS}}
\eeq
with effective absolute constants $c_2, c_3 > 0$.

We remark that \eqref{eq:2.1} and the second inequality of \eqref{eq:2.2} imply
\beq
\label{eq:2.3}
r \geqslant c_4 \log^2 N
\eeq
and for the greatest prime factor $q_0$ of $r$
\beq
\label{eq:2.4}
q_0 \geqslant 2 \log_2 N - c_5 > \log_2 N \ \text{ if } \ N > N_0,
\eeq
with effective absolute constants $c_4, c_5 > 0$.

We will slightly reformulate Theorem~4.2 of \cite{BFM} which itself is an improved reformulation of the original Maynard--Tao theorem.
We remark that in order to obtain Theorems \ref{th:1}--\ref{th:3} (with a constant $C$ larger than $50$, respectively with a proportion less than $1/50 = 2\%$) one could use also the more complicated method of Zhang \cite{Zha}.
However, to obtain a new proof of Theorem~\ref{th:4} we need the Maynard--Tao method.
Also the result stated below uses clearly the Maynard--Tao method.

Let $P^+(n)$ denote the largest prime factor of~$n$.

\begin{theo}[Maynard--Tao]
Let $k = k_m$ be an integer, $\ve = \ve(k,n) > 0$ be sufficiently small, $N > N_0(\ve, k, m)$.
Further, let
\beq
\label{eq:2.5}
k + 1 < C_6(\ve) < h_1 < h_2 < \dots < h_k \leqslant N, \ \mathcal H = \mathcal H_k = \{h_i\}_{i = 1}^k \text{ admissible},
\eeq
\beq
\label{eq:2.6}
\Delta(\mathcal H) := \prod_{1 \leqslant i < j \leqslant k} (h_j - h_i), \ \ \left(q_0 \prod_{i=1}^k h_i , \Delta(\mathcal H)\right) = 1, \ \ \Delta(\mathcal H) < N^\ve.
\eeq
\beq
\label{eq:2.7}
\text{For } \ m = 2 \ \text{ let } \ k = 50
\eeq
and generally let
\beq
\label{eq:2.8}
k_m = C_7 e^{5m}
\eeq
with suitably chosen constants $C_7$ and $C_6(\ve)$, depending on~$\ve$.
Then we have at least $m$ primes among $n + \mathcal H_k = \{n + h_i\}_{i = 1}^k$ for some $n \in (N, 2N]$.
\end{theo}

\begin{remark}
\label{rem:1}
In the proof of our Theorems~\ref{th:1}--\ref{th:4} we will have
\beq
\label{eq:2.9}
h_k \leqslant g(N) \log N < \log^2 N,
\eeq
thus the second condition of \eqref{eq:2.6} will be trivially fulfilled.
\end{remark}

\begin{remark}
\label{rem:2}
In the mentioned applications we will choose the values $h_i$ as primes, so the first condition of \eqref{eq:2.6} will be equivalent to
\beq
\label{eq:2.10}
q_0 \nmid h_j - h_i, \ h_t \nmid h_j - h_i \ \text{ for any } \ t \in [1,k], \ 1 \leqslant i < j \leqslant k.
\eeq
Since in the applications the only other condition will be with some functions $\xi_i(N)$ to have
\beq
\label{eq:2.11}
 h_i = (1 + o(1))\xi_i(N), \ \ \ \xi_i(N) \ll g(N) \log N
\eeq
it will make no problem to choose step by step primes $h_i$ satisfying \eqref{eq:2.9}--\eqref{eq:2.10}.
Also $h_i \in \mathcal P$, $h_i > k$ assures that $\mathcal H_k$ is admissible.
\end{remark}

The Maynard--Tao method assures the existence of at least $m$ primes among numbers of the form
\beq
\label{eq:2.11masodik}
n + h_i \ (1 \leqslant i \leqslant k) \ \text{ with } \ n \equiv z \ (\text{\rm mod } W)
\eeq
with any $z \in [1, W]$ and for some $n \in (N, 2N]$, $N$ sufficiently large, if
\beq
\label{eq:2.12}
\left(\prod_{i = 1}^k (z + h_i), W\right) = 1.
\eeq
The pure existence of such a $z$ follows from the admissibility of $\mathcal H_m$ but its actual choice is crucial in the applications.

In order for the method of Maynard--Tao and Banks--Freiberg--Maynard \cite{BFM} to work we must assure still (see \cite{BFM}) with a sufficiently large $C_8(\ve)$
\beq
\label{eq:2.13}
\Delta(\mathcal H) = \prod_{1\leqslant i < j \leqslant k} (h_j - h_i)\, \Big| \, W; \ \ \ \prod_{p \leqslant C_8(\ve)} p \, \Big |\, W
\eeq
and for the possibly existing greatest prime factor $q_0$ of the possibly existing exceptional modulus $r$,
\beq
\label{eq:2.14}
q_0 \nmid W \ \ (\text{if such a modulus $r$, and so $q_0$ exists});
\eeq
further ($P^+(n)$ will denote the greatest prime factor of $n$)
\beq
\label{eq:2.15}
P^+(W) < N^{\ve / \log_2 N}, \ W < N^{2\ve}.
\eeq
In the applications \eqref{eq:2.9} will assure $q_0 \nmid \Delta(\mathcal H) \prod\limits_{p \leqslant C_8(\ve)} p$.

If we succeed to show the existence of a pair $(z, W)$ with \eqref{eq:2.12}--\eqref{eq:2.15} and the crucial additional property that with a suitable $c_9(\ve)$
\beq
\label{eq:2.16}
(z + s, W) > 1 \ \text{ if } \ s \notin \mathcal H, \ \ 1 < s \leqslant c_9(\ve) g(N) \log N,
\eeq
then we can assure that all numbers $z + s$ with \eqref{eq:2.16} have a prime divisor $p \mid W$.
Consequently all $n + s$, $s \neq h_i$, $s \in (1, c_9(\ve)g(N) \log N]$ will be composite if $n \in (N, 2N]$.

In order to achieve this, we will use the Erd\H{o}s--Rankin method.
After this we can show Theorems \ref{th:1}--\ref{th:4} with suitable choices of $\mathcal H_k$.

We will choose the following parameters ($p$ will always denote primes), $\mathcal H = \mathcal H_i$, $c_{10}(\ve) = 2 c_9(\ve)/\ve$,
\beq
\label{eq:2.17}
\mathcal L = \ve \log N, \ \ v = \log^3 \mathcal L, \ \ U = c_{10}(\ve)g(e^{\mathcal L})\mathcal L > c_9(\ve) g(N) \log N,
\eeq
\beq
\label{eq:2.18}
y = \exp \left(\frac1{k + 5} \log \mathcal L \log_3 \mathcal L / \log_2 \mathcal L\right), 
\eeq
\beq
\label{eq:2.19}
P_1 = \underset{p \leqslant v}{\prod\nolimits^*} p,
\eeq
\beq
\label{eq:2.20}
P_2 = \underset{v < p \leqslant y}{\prod\nolimits^*} p,
\eeq
\beq
\label{eq:2.21}
P_3 = \underset{y < p \leqslant \mathcal L/2 }{\prod\nolimits^*} p,
\eeq
\beq
\label{eq:2.22}
P_4 = \underset{\mathcal L/2 < p \leqslant \mathcal L}{\prod\nolimits^*} p,
\eeq
where $\underset{p}{\prod^*}$ means
\beq
\label{eq:2.23}
p \not\in \mathcal H':= \mathcal H\cup \{q_0\}.
\eeq
Further, let
\beq
\label{eq:2.24}
W_0 = P_1 P_2 P_3 P_4, \ \ W = \biggl[\underset{p \leqslant \mathcal L}{\prod\nolimits^*} p , \Delta_0(\mathcal H)\biggr] = \bigl[P_1 P_2 P_3 P_4, \Delta_0 (\mathcal H)\bigr]
\eeq
where $\Delta_0(\mathcal H)$ denotes the squarefree part of $\Delta(\mathcal H)$.
This choice of $W$ clearly satisfies both conditions of \eqref{eq:2.15} by the Prime Number Theorem if we additionally require the condition valid in all applications:
\beq
\label{eq:2.25}
h_k \leqslant \log^2 N.
\eeq

\section{The application of the Erd\H{o}s--Rankin method}
\label{sec:uj4}

We will choose the congruence class $z$ modulo any prime divisor of $W$,
which finally determines $z$ $\text{\rm mod } W$.
Let $p$ denote always primes; further let us choose
\beq
\label{eq:3.1}
z \equiv 0 \ (\text{\rm mod } P_1 P_3).
\eeq
This implies by $v \mathcal L / 2 > U$ that
\beq
\label{eq:3.2}
(z + s, P_1 P_3) = 1 \qquad (1 < s \leqslant U)
\eeq
if and only if $(s, P_1 P_3) = 1$, that is, if and only if either
\beq
\label{eq:3.3}
s = pq_0^\alpha \prod_{i = 1}^k h_i^{\alpha_i} \ \ (\alpha \geqslant 0, \ \alpha_1 \geqslant 0, \dots, \alpha_m \geqslant 0) \ \text{ and } \ p > \mathcal L/2
\eeq
or
\beq
\label{eq:3.4}
s \text{ is composed only of primes } \ p \,\Big |\, P_2 q_0 \prod_{i = 1}^k h_i.
\eeq

The first step is to estimate the number $A_0$ of numbers $s$ satisfying \eqref{eq:3.4}.
This will be relatively easy since

\begin{itemize}
\item[(i)]
we have an upper estimate for $y$-smooth numbers by the results of Dickman and (in a refined form) of de Bruijn \cite{Bru}.
We quote a suitable result of \cite{Bru} as our Lemma~\ref{lem:1} in a simpler form as given in \cite{Mai1};
\item[(ii)]
the additional factor $q_0^\alpha \prod\limits_{i = 1}^k h_i^{\alpha_i}$ leaves the asymptotic for numbers of the form \eqref{eq:3.3} and \eqref{eq:3.4} below $U$ nearly unchanged.
\end{itemize}

\begin{lemma}
\label{lem:1}
Let $\Psi(x,y)$ denote the number of positive integers $n \leqslant x$ which are composed only of primes $\leqslant y$.
For $y \leqslant x$, $y \to \infty$, $x \to \infty$ we have
\beq
\label{eq:3.5}
\Psi(x,y) \leqslant x \exp \left[ - \log x \frac{\log_3 y}{\log y} + (1 + o(1)) \log_2 y\right].
\eeq
\end{lemma}

This is a slightly simplified form of Lemma~5 of \cite{Mai1}.

Applying \eqref{eq:3.5} and taking into account that $\alpha, \alpha_i \ll \log U \sim \log \mathcal L$, we obtain by the choice of $y$ in \eqref{eq:2.18} for the number of $s \leqslant U$ with \eqref{eq:3.4} the upper estimate
\begin{gather}
\label{eq:3.6}
(1 + o(1))(\log \mathcal L)^{k + 1} U \exp \left[- \frac{(k\! +\! 5\! +\! o(1)) \log\mathcal L \log_3 \mathcal L}{\log \mathcal L \log_3 \mathcal L / \log_2 \mathcal L} + (1 + o(1)) \log_2 \mathcal L \right]\! \ll\\
\ll \frac{U}{\log^2 U} \leqslant \frac{C_{11}\bigl(\pi(\mathcal L) - \pi(\mathcal L/2)\bigr)}{\log U}
\nonumber
\end{gather}
which will be negligible compared with the numbers of $s \leqslant U$ with \eqref{eq:3.3}.
This means that integers with \eqref{eq:3.4} can be later sieved out by a tiny portion of primes dividing $P_4$ in \eqref{eq:2.22}.
On the other hand, the number of integers with \eqref{eq:3.3} is much larger than $\pi(\mathcal L)$, although just slightly larger than
\hbox{$\pi(U) - \pi(\mathcal L/2)$}, which corresponds to the case $\alpha = \alpha_1 = \dots = \alpha_k = 0$ in \eqref{eq:3.3}.
We have, namely, by the Prime Number Theorem,
\begin{align}
\label{eq:3.7}
A_0' &= \sum_{\alpha, \alpha_1, \dots, \alpha_k \geqslant 0} \pi \left(\frac{U}{q_0^\alpha \prod\limits_{i = 1}^k h_i^{\alpha_i}}\right) - \pi \left(\frac{\mathcal L/2}{q_0^\alpha \prod\limits_{i = 1}^k h_i^{\alpha_i}}\right) \\
&\leqslant \frac{(1 + o(1)) U}{\log U} \left(1 + \frac1{q_0} + \frac1{q_0^2} + \dots\right) \prod_{i = 1}^k \left(1 + \frac1{h_i} + \frac1{h_i^2} + \dots \right)\nonumber \\
&\sim \frac{U}{\log U} \left(1 - \frac1{q_0}\right)^{-1} \prod_{i = 1}^m \left(1 - \frac1{h_i} \right)^{-1} \leqslant \frac{2U}{\log U}
\nonumber
\end{align}
if $C_6(\ve)$ was chosen sufficiently large depending on~$k$.

We will choose the residue class $z_{p_j}$ $(\text{\rm mod }p_j)$ for all $p_j \mid P_2$ consecutively for all primes.
We have to take care in the $j$\textsuperscript{th} step that
\beq
\label{eq:3.8}
z_{p_j} + h_i \not \equiv 0 \ \ (\text{\rm mod }p_j)
\eeq
should hold; further, the additional property that at the $j$\textsuperscript{th}
step we choose the residue class $z_p$ $\text{\rm mod }p$ so that it should sieve out the maximal number of remaining elements
from the remaining set of $s$'s of cardinality $A_{j - 1}'$.
We distinguish two cases:

(i) if before the $j$\textsuperscript{th} step the number of $s$'s satisfies
\beq
\label{eq:3.9}
A_{j - 1}' \leqslant \frac{\mathcal L}{5 \log \mathcal L} \ \ \left( < \frac{\pi(\mathcal L) - \pi(\mathcal L/2)}{2}\right),
\eeq
then we stop the choice of new $z_p$'s.

Otherwise, if \eqref{eq:3.9} is false, then we have in total at most $(\log_2 \mathcal L)^{k + 1}$ possibilities for $\alpha$, $\{\alpha_i\}_{i = 1}^k$, since $\dfrac{U}{p_j} = o(\log \mathcal L)$ and even neglecting the primality of $p$ we have in total at most
\beq
\label{eq:3.10}
k \left\lceil \frac{U}{y}\right\rceil (\log_2 \mathcal L)^{k + 1} < \frac{\mathcal L}{\log^{10} \mathcal L} < \frac{A_{j - 1}'}{2 (\log \mathcal L)^8}
\eeq
numbers $s$ in forbidden residue classes \eqref{eq:3.8}.
This means that choosing the residue class $z_{p_j}$ so that we avoid the $k$ forbidden residue classes
$h_1, \dots, h_k$ but sieve out afterwards as many elements as possible, we obtain after the next step
\begin{align}
\label{eq:3.11}
A_j' &< A_{j - 1}' - \frac{A_{j - 1}' \left(1 - \frac1{2(\log \mathcal L)^8}\right)}{p_j - k} \\
&< A_{j - 1}' \left(1 - \frac{1 - \frac1{2(\log\mathcal L)^8}}{p_j}\right)\nonumber\\
&< A_j \left(1 - \frac1{p_j}\right)^{1 - (\log\mathcal L)^{-8}}.\nonumber
\end{align}
By Mertens' theorem, \eqref{eq:2.17}--\eqref{eq:2.18} and \eqref{eq:3.7}, we obtain a final residual set
 (after at most $\pi(y) - \pi(v)$ steps of $s$'s) of size at most
\begin{align}
\label{eq:3.12}
A^* &< A_0'\prod_{v < p \leqslant y} \left(1 - \frac1{p}\right)^{1 - (\log \mathcal L)^{-8}} \sim A_0' \left(\frac{\log v}{\log y} \right)^{1 - (\log \mathcal L)^{-8}} \\
&\sim A_0' \left(\frac{3(k + 5) \log_2^2 \mathcal L}{\log \mathcal L \log_3 \mathcal L}\right) < \frac{7kU}{\log U g(e^{\mathcal L})} < \frac{7kc_{10}(\ve)\mathcal L}{\log \mathcal L}
\nonumber\\
&< \frac{\pi(\mathcal L) - \pi(\mathcal L/2)}{3}.
\nonumber
\end{align}

This means that taking into account that the total number of $s$'s with \eqref{eq:3.4} is by \eqref{eq:3.6} a negligible portion of the above remaining quantity (even without sieving them out by the above procedure), we obtain finally that using the prime factors of $P_2$, with a suitable choice of $z_p$ for these primes we can already reach for the numbers $s$ with $1 < s \leqslant U$, $(z + s, P_1 P_3) = 1$ apart from an exceptional set $S$ of size at most $\bigl(\pi(\mathcal L) - \pi(\mathcal L / 2)\bigr)/2$ the crucial relations
\beq
\label{eq:3.13}
(z + s, P_2) > 1 \ \text{ if } \ s \notin \mathcal H, \ \ s \notin S
\eeq
and
\beq
\label{eq:3.14}
(z + s, P_2) = 1 \ \text{ if } \ s \in \mathcal H, \ \ {s \notin S}.
\eeq
Let
\beq
\label{eq:3.15}
S' = S \setminus \mathcal H.
\eeq
Then by $|S| < \bigl(\pi(\mathcal L) - \pi(\mathcal L / 2)\bigr)/2$ we can easily find for any $s \in S$ a suitable prime $p \mid P_4 = \bigl\{ p \in (\mathcal L/2, \mathcal L], p \notin \mathcal H\bigr\}$ with $p \nmid \prod\limits_{i = 1}^k (s - h_i)$ and consequently a $z_p(\text{\rm mod } p)$ with
\beq
\label{eq:3.16}
z_p + s \equiv 0 \ (\text{\rm mod }p) \text{ for } s \in S', \ \ z_p + t \not \equiv 0\ (\text{\rm mod }p) \text{ for } t \in \mathcal H.
\eeq

Thus we need still to determine $z$ $\text{\rm mod }p$ for those primes which were not used before.
These primes belong to one of the following categories (cf.\ \eqref{eq:2.5}):

\begin{itemize}
\item[(i)] \ $p \in \mathcal H \longrightarrow p > C_6(\ve)$,

\item[(ii)] \ $p = q_0 \longrightarrow p \geqslant \log_2 N $ \ (if $q_0 \leqslant W$),

\item[(iii)] \ the remaining parts of unused $p \mid P_4 \longrightarrow p > \mathcal L/2$,

\item[(iv)] \ $p \,\Big| \, \dfrac{W}{W_0} \longrightarrow p > \mathcal L$.
\end{itemize}

Since we have for $1 < s \leqslant U$, $s \notin \mathcal H$ already by the earlier choices
\beq
\label{eq:3.17}
(z + s, P_1 P_2 P_3 P_4) > 1
\eeq
this property will be valid independently from the further choices of $z_p$ and so the condition $z + s$ composite for $1 < s \leqslant U$, $s \notin \mathcal H$ will be true at the end as well.

So we have only to assure that for the primes in (i)--(iv) we should have
\beq
\label{eq:3.18}
z_p \not \equiv -t \ (\text{\rm mod }p) \ \text{\rm if } \  t \in \mathcal H.
\eeq
But this makes no problem since $|\mathcal H| = k < \min \bigl\{C_6(\ve), \log_2 N, \mathcal L/2\bigr\}$ by \eqref{eq:2.5}.

So we finally determined a $z$ $\text{\rm mod }W$ with the property that for $s \in (1, U]$ we have
\beq
\label{eq:3.19}
(z + s, W) = 1 \ \text{ if and only if } \  s \in \mathcal H.
\eeq
Consequently if $n \equiv z$ $(\text{\rm mod }W)$, then
\beq
\label{eq:3.20}
(n + s, W) = 1 \ \text{ if and only if } \  s \in \mathcal H.
\eeq

\section{Proofs of Theorems \ref{th:1}--\ref{th:4}}
\label{sec:uj5}

We summarize the results of Sections~\ref{sec:uj3} and \ref{sec:uj4} with the aim of applications in Theorems \ref{th:1}--\ref{th:4}.

Let $k, m, \ve$ be chosen satisfying \eqref{eq:4.1}, let $\ve = \ve(k, m) > 0$ be a sufficiently small constant, $N > N_0(\ve, k, m)$:
\beq
\label{eq:4.1}
k = 50, \ \ m = 2 \ \text{ or } \ k_m = C_7 e^{5m}.
\eeq
Let $\mathcal H = \{h_i\}_{i = 1}^k \subset \mathcal P$ satisfying
\beq
\label{eq:4.2}
\bigl[h_t, q_0\bigr] \nmid \prod_{1 \leqslant i < j \leqslant k} (h_j - h_i) \ \text{ for any } \ t \in [1, m], \ q_0 \text{ defined in \eqref{eq:2.1}--\eqref{eq:2.4},}
\eeq
\beq
\label{eq:4.3}
k + 1 < C_6(\ve) < h_1 < \dots < h_k \leqslant \log^2 N,
\eeq
\beq
\label{eq:4.4}
U = c_9(\ve)g(N) \log N, \ \ g(N) = \frac{\log_2 N \log_4 N}{\log_3^2 N}, \ N' = \pi(N).
\eeq
Then we can find suitable values of
\beq
\label{eq:4.5}
W < N^\ve, \ z(\text{\rm mod } W)
\eeq
and an $n \in [N', 2N')$ such that we have at least $m$ primes among $n + h_i$ and all numbers of the form $n + s$ are
composite if $s \in (1, U] \setminus \mathcal H$.

\begin{rema}
This implies that all (at least $m$) primes in the interval $(n + 1, n + U)$ are of the form $n + h_j$, $h_j \in \mathcal H$.
\end{rema}

\begin{rema}
We used the introduction of the variable $\ve$ since it was formulated in this way in \cite{BFM}.
However, since there is an exact connection \eqref{eq:4.1} between $k$ and $m$ and $\ve$ depends just on $k$ and $m$, in the applications we can write in \eqref{eq:4.3} $C_6'(k)$ instead of $C_6(\ve)$, $C_8'(k)$ in place of $C_8(\ve)$ in \eqref{eq:2.13} and $c_9'(k)$ instead of $c_9(\ve)$ in \eqref{eq:4.4}, $c_{10}'(k)$ instead of $c_{10}(\ve)$ before \eqref{eq:2.17}, further $c_{11}'(k)$ in place of $\ve$ in \eqref{eq:4.5}.
Similarly we can choose $\mathcal L = c_{12}(k) \log N$ with a small $c_{12}(k)$ in \eqref{eq:2.17}.
Additionally, if $m = 2$, $k = 50$, these are just absolute constants (which is the case in Theorems \ref{th:1}, \ref{th:2} and \ref{th:3}).
The elimination of $\ve$ in this part of the proof will also increase clarity since the condition \eqref{eq:1.5} for the function $f(n)$ contains a parameter $\ve$ too.
\end{rema}

In order to prove Theorem~\ref{th:1}, suppose, in contrary to its assertion, that we have a sequence of 50 positive numbers $c_\nu^*$, $\delta_\nu$ $(1\leqslant \nu \leqslant 50)$ satisfying with two constants $c^*$, $N^* > 0$
\beq
\label{eq:4.6}
J_\nu := \bigl[c_\nu^*, c_\nu^* + \delta_\nu\bigr], \ \ c_\nu^* > 4 \delta_\nu > 20 c_{\nu + 1}, \ \ c_1^* < c^*,
\eeq
\beq
\label{eq:4.7}
\left\{\frac{d_n}{f(n)}\right\}_{n = N^*}^\infty \cap \biggl(\bigcup_{\nu = 1}^{50} J_\nu \biggr) = \emptyset.
\eeq

Let
\beq
\label{eq:4.8}
I_\nu(n) := \Bigl[c_\nu^* f(n), \bigl(c_\nu^* + \delta_\nu\bigr) f(n)\Bigr] \ \text{ for } \nu = 1,2, \dots, 50.
\eeq
Then
\beq
\label{eq:4.9}
d_n \notin \bigcup_{\nu = 1}^{50} I_\nu(n) \ \text{ for } \ \nu = 1,2, \dots, 50, \ n \in [N', 2N'), \ N' > N^*.
\eeq
We will choose now the primes $h_1 < h_2 < \dots < h_{50}$ consecutively, satisfying \eqref{eq:4.2}--\eqref{eq:4.3} and a sufficiently small $\ve > 0$, $N' > \max \bigl(N_0(\ve), N^*\bigr)$
\beq
\label{eq:4.10}
h_\nu \in I_{51 - \nu}' (n) := \left[\left(c_{51 - \nu}^* + \frac{\delta_{51 - \nu}^*}{2} \right) (1 + \ve) f(N'),
\bigl(c_{51 - \nu}^* + \delta_{51 - \nu}\bigr) (1 - \ve) f(N')\right].
\eeq

This choice (for $\nu = 1, \dots, 50$) is easily assured by the Prime Number Theorem if $\ve$ was chosen sufficiently small, $N_0(\ve)$ sufficiently large depending on all $c_\nu^*$, $\delta_\nu$ ($1 \leqslant \nu \leqslant 50$) and $\ve$.
So we have for $1 \leqslant \nu < \mu \leqslant 50$ for large enough $N'$
\begin{align}
\label{eq:4.11}
h_\mu - h_\nu \in \, &\biggl[\left( c_{51 - \mu}^* + \frac{\delta_{51 - \mu}}{2} - 2 c_{52 - \mu}^*\right) (1 + \ve) f(N'),
\\
&\bigl(c_{51 - \mu}^* + \delta_{51 - \mu}\bigr)(1 - \ve) f(N')\biggr] : = I_{51 - \mu}^* (n) \subset I_{51 - \mu}(n).
\nonumber
\end{align}
This contradicts to \eqref{eq:4.9} since we have for at least one pair of consecutive primes
\beq
\label{eq:4.12}
d_n = h_\mu - h_\nu, \ \ n \in [N', 2N') \ \ \text{q.e.d.}
\eeq

Now we turn to the proof of Theorem~\ref{th:2}.
Let us suppose that we have $50$ functions with $f_i(x) \nearrow \infty$, satisfying \eqref{eq:1.5}, \eqref{eq:1.9}, \eqref{eq:1.11} and $50$ intervals
\beq
\label{eq:4.13}
J_\nu := \bigl[c_\nu^*, c_\nu^* + \delta_\nu\bigr], \ \ \ I_\nu(n) := \bigl[c_\nu^* f_\nu(n), (c_\nu^* + \delta) f_\nu(n)\bigr]
\eeq
such that with a sufficiently large $N^*$ we have
\beq
\label{eq:4.14}
d_n \notin \bigcup_{\nu = 1}^{50} I_\nu(n) \ \text{ for } \ \nu = 1, \dots, 50, \ \ n \in [N', 2N'), \ N' > N^*.
\eeq

Then, analogously to \eqref{eq:4.10}--\eqref{eq:4.11} we can choose the primes $h_1 < h_2 < h_{50}$ with \eqref{eq:4.2}--\eqref{eq:4.3},
$\ve > 0$ and $N > \max\bigl(N_0(\ve), N^*\bigr)$ so that
\beq
\label{eq:4.15}
h_\nu \in I_\nu'(n) := \left[ \left(c_\nu^* + \frac{\delta_\nu}{2}\right)(1 + \ve) f_\nu(N'), (c_\nu^* + \delta_\nu) (1 - \ve) f_\nu(N')\right].
\eeq
This implies for sufficiently large $N'$ by \eqref{eq:1.11} for any $1 \leqslant \nu < \mu \leqslant 50$
\beq
\label{eq:4.16}
h_\mu - h_\nu \in \left[c_\mu^*(1 + \ve) f_\mu(N), \bigl(c_\mu^* + \delta_\mu\bigr)(1 - \ve) f_\mu(N)\right] \subset  I_\mu(n)
\eeq
and we obtain again a contradiction to \eqref{eq:4.14}. \ q.e.d.

The first assertion of Theorem~\ref{th:3} follows immediately from our summary in \eqref{eq:4.1}--\eqref{eq:4.5} if we choose simply
$h_1 < h_2 < \dots < h_{50}$ satisfying \eqref{eq:4.2}--\eqref{eq:4.3} and for $n \in [N', 2N')$ with
\beq
\label{eq:4.17}
h_i = \beta_i f(N')\left(1 + O\left(\frac1{\log_2 N'}\right)\right).
\eeq
The consequence about the at least 2\% density of $J_f$ follows in the same way as in the proof of Corollary~1.2 of \cite{BFM}.

Finally, Theorem~\ref{th:4} is also an obvious corollary of our summary \eqref{eq:4.1}--\eqref{eq:4.5}.
Namely, if for $n \in [N', 2N')$ we choose for $i \in [1, k_m]$
\beq
\label{eq:4.18}
h_i = i \cdot \frac{U}{k_m + 1} \left(1 + O\left(\frac1{\log U}\right)\right) \ \ \ \bigl(k_m = C_7 e^{5m}\bigr)
\eeq
with \eqref{eq:4.2}--\eqref{eq:4.5} then we have in the interval $(n + 1, n + U)$ at least $m$ and at most $k_m$ primes, all among $n + h_i$ $(1 \leqslant  i \leqslant k_m)$.
Consequently we get at least $m$ consecutive primegaps each of size at least
\beq
\label{eq:4.19}
\frac{U}{2(k_m + 1)} \geqslant c_{12}(k) g(n) \log n \ \text{ with } \ g(n) = \frac{\log_2 n \log_4 n}{\log_3^2 n}. \ \ \text{ q.e.d.}
\eeq

\section{Proofs of Theorems 1'--4' and Corollary 1'}
\label{sec:uj6}

We first remark that apart from Theorem~\ref{th:2} we worked in all proofs (cf.\ our present Section~\ref{sec:uj5}) within a given interval $[N', 2N']$ where $N'$ was any sufficiently large constant and we worked with an $\mathcal H_k$ tuple satisfying
\beq
\label{eq:6.1}
\mathcal H_k = \{h_i\}_{i = 1}^k \ \text{ with } \ h_i \asymp f(N').
\eeq

Thus we will consider first the proofs of Theorems 1', 3', 4'.
We will distinguish two cases as follows.

\smallskip
\emph{Case 1}.
\[
f(N') < \log N'(\log_2 N')^{1/2}.
\]

In this case the assertions of Theorems 1', 3', 4' follow directly from Theorems~\ref{th:1}, \ref{th:3}, \ref{th:4} for the specific interval $[N', 2N']$.

\smallskip
\emph{Case 2}.
\[
f(N') \geqslant  \log N'(\log_2 N')^{1/2}.
\]

In this case we will use the method of \cite{May2}, and will describe the needed changes compared to \cite{May2}.
We will use \eqref{eq:6.1} which in this case implies
\beq
\label{eq:6.2}
h_i \gg (\log N') (\log_2 N')^{1/2} \ \ \ (i = 1,2, \dots, k).
\eeq

In order to follow \cite{May2} we will change our notation and choose with a given small $\varepsilon_0$
\beq
\aligned
 z &= \varepsilon_0 \log N', \ \ x = \mathcal L, \ \ P_y = \prod_{p \leqslant y} p, \\
y &= \exp \left( \frac{(1 - \varepsilon_0) \log x \log_3 x}{\log_2 x} \right), \ \ z = \frac{x}{\log_2 x}, \ \ U = C_U \frac{x \log y}{\log_2 x},
\endaligned
\label{eq:6.3}
\eeq
where $C_0$ is an arbitrarily large constant as in \cite{May2}, independent of $\varepsilon_0$.

In contrast to Section~\ref{sec:uj2} and in accordance with \cite{May2} we will choose the residue classes $a_p(\text{\rm mod }p)$, in the first step for $p \leqslant z$, $p \neq q_0$ (the greatest prime factor of the eventually existing single exceptional modulus, as in Sections \ref{sec:uj3}--\ref{sec:uj4})
\beq
\label{eq:6.4}
a_p = 0 \ \text{ for every prime } \ p \in (y, z], \ \ p\neq q_0,
\eeq
\beq
\label{eq:6.5}
a_p = 1 \ \text{ for every prime } \ p \leqslant y, \ \ p \neq q_0.
\eeq

After removing elements of $[1, U]$ in these residue classes we obtain the set $\mathcal R \cup \mathcal R' \cup \widetilde{\mathcal R} \cup \widetilde{\mathcal R}'$, where
\begin{align}
\label{eq:6.6}
\mathcal R &= \left\{ m p \leqslant U: \ p > z, \ m \text{ is $y$-smooth, } (mp - 1, P_y) = 1 \right\}, \\
\widetilde{\mathcal R} &= \left\{ mpq_0 \leqslant  U : \ p > z, \ m \text{ is $y$-smooth, }
(mpq_0 - 1, P_y) = 1 \right\}, \nonumber\\
\mathcal R' &= \left\{m \leqslant U:\ m \text{ is $y$-smooth, } (m - 1, P_y) = 1 \right\}, \nonumber\\
\widetilde{\mathcal R}' &= \left\{mq_0 \leqslant U: \ m \text{ is $y$-smooth, } (mq_0 - 1, P_y) = 1 \right\}. \nonumber
\end{align}

We obtain by Lemma~\ref{lem:1} of Section~\ref{sec:uj4} similarly to Lemma~2 of \cite{May2}
\beq
\label{eq:6.7}
\left| \mathcal R'\cup \widetilde{\mathcal R}'\right| \ll \frac{x}{(\log x)^{1 + \varepsilon}} .
\eeq

Again similarly to Lemma~3 of \cite{May2} we have now for $V \in \left[z + z / \log x, x(\log x)^2\right]$
\beq
\label{eq:6.8}
\#\left\{z < p \leqslant V : (mp - 1, P_y) = 1 \right\} = \frac{V - z}{\log x} \prod_{\substack{p \leqslant y\\ p \nmid m}} \frac{p - 2}{p - 1} (1 + o(1))
\eeq
and in particular for even $m \leqslant U(1 - 1 / \log x)/z$
\beq
\label{eq:6.9}
\left|\mathcal R_m\right| = \frac{2 e^{-\gamma}  U (1 + o(1))}{m(\log x)(\log y)} \left(\prod_{p > 2} \frac{p(p - 2)}{(p - 1)^2}\right)
\left( \prod_{p \mid m, p > 2} \frac{p - 1}{p - 2}\right).
\eeq

Further the same methods show that by $q_0 \geqslant  \log_2 N'$
\beq
\label{eq:6.10}
\left|\widetilde{\mathcal R}_m \right| \ll \frac{U}{q_0 m \log x \log y} \ll \frac{U}{m \log^2 x \log y},
\eeq
which is negligible compared to \eqref{eq:6.9}.
By $\mathcal R_m$ and $\widetilde{\mathcal R}_m$, resp., we denoted the terms of $\mathcal R$ and $\widetilde{\mathcal R}$, resp., which contain the specific parameter $m$ in \eqref{eq:6.6} as in (2.8) of \cite{May2}.

After this initial choice of $a_p$ for $p \leqslant z$ we try to choose $a_p$ modulo $p$ for all $p \leqslant x$, 
$p \neq q_0$
 in such a way that for the common solution $a$ modulo $W$ of the congruences $a \equiv a_p \ (\text{\rm mod }p)$ for
\beq
\label{eq:6.11}
p \mid W := \frac{P_x}{q_0}
\eeq
we should have, similarly to \eqref{eq:2.17}
\beq
\label{eq:6.12}
(a + s, W) > 1 \ \text{ if } \ s \notin \mathcal H, \ \ 1 < s \leqslant U.
\eeq

The choice of $a_p$ for $z < p \leqslant x$ will follow closely that of \cite{May2} with the following changes.
We will choose in the applications (see Section~\ref{sec:uj5}) our $\mathcal H$ with
\beq
\label{eq:6.13}
h_i \in \mathcal P, \ \ (h_i - 1, q_0 P_y) = 1
\eeq
and use \eqref{eq:6.2} additionally.
This means that $h_i \equiv 0 \ (\text{\rm mod }p)$ for some $p \in (y, z]$ will not occur and we will have
\beq
\label{eq:6.14}
h_i \in \mathcal R_1 \ \ (i = 1,2, \dots, k).
\eeq
(We note that $\mathcal R_m \cap \mathcal R_{m'} = \emptyset$ if $m \neq m'$, so $h_i \notin \mathcal R_m$ for $m > 1$.)

Since any other essential requirements for $h_i$ in Sections \ref{sec:uj3}--\ref{sec:uj5} are concerning only the size of $h_i$ requiring
\beq
\label{eq:6.15}
\xi_i(N) \leqslant h_i \leqslant \xi_i (N) (1 + \eta_i)
\eeq
for some sufficiently small $\eta_i$ independent of $N$, with some functions $\xi_i(N)$, we can always fulfil these conditions with proper choice of the primes $h_i$ satisfying \eqref{eq:6.13}, due to the relations \eqref{eq:6.8}--\eqref{eq:6.9}, which mean that we have sufficiently large sets to choose $\{h_i\}_{i = 1}^k $.

After choosing our set $\mathcal H = \{h_i\}_{i = 1}^k$ satisfying the requirements of Sections \ref{sec:uj3}--\ref{sec:uj5} for a given value of $N$ we will denote it by $\mathcal H_k^* = \{h_i^*\}_{i = 1}^{k^*}$ and consider it fixed.
After this we will choose $a_p$ for $p \mid W$, $p > z$ in a somewhat different way from that of \cite{May2}.
The difference affects only the case $m = 1$ and can be described as follows.
We will choose the set $\mathcal H = \{h_1, \dots, h_k\}$ in \cite{May2} disjoint to our $\mathcal H_{k^*}^*$ (and also $k$ will be sufficiently large compared to $k^*$ while our $k^*$ will be equal to $50$, $99$ or $k = k(m)$ in Theorem~4').
The change in the choice of the probabilities of choosing $a$ $\text{\rm mod } q \in I_1 \subseteq [x/2, x]$ will be that in contrast to (4.1) of \cite{May2} we will set ($\mu^*_{1,q}(a)$ will denote the new probabilities, $\alpha_{1,q}^* $
the new normalizing number to have $\sum\limits_{a(q)} \mu_{1,q}^* (a) = 1$)
\begin{align}
\label{eq:6.16}
\mu_{1,q}^*(a) &:= 0 \ \text{ if } \ \exists i \in [1, k^*], \ \ a + h_i^* \equiv 0 \ (\text{\rm mod }q),\\
\label{eq:6.17}
\mu_{1,q}^*(a) &:= \mu_{1,q}(a) \cdot \frac{\alpha_{1,q}^*}{\alpha_{1,q}} \ \ \text{ otherwise}.
\end{align}

In this way we can avoid that by the random choice of $a_q$ modulo $q$ we should have
$q \mid n + h_i^*$ for $n \equiv a \ (\text{\rm mod } W)$ for some $i \in [1, k^*]$, since we give $0$ probability to those $a_q$ in \eqref{eq:6.16}.
Naturally we have to rescale the remaining probabilities as done in \eqref{eq:6.17}, which actually slightly increases all remaining probabilities.
This means that none of the $h_i^* \in \mathcal R_1$ will be sieved out (with probability $1$) by the above random sieve procedure.
On the other hand we have to show that, similarly to the end of Section~6 on p.~13 of \cite{May2}, for all but $o_k\left(\left|\mathcal R_1 \right|\right)$ primes $p_0 \in \mathcal R_1$ the expected number $\sum\limits_q \mu_{1,q}^*(p_0)$ of times $p_0 \in \mathcal R_1$ is chosen will remain $\gg \delta \log k$ if $p_0 \notin \mathcal H^*$ as in \cite{May2} in case of the original choice of $\mu_{1,q}(p_0)$.

If for a given $q$
\beq
\label{eq:6.18}
p_0 \not\equiv h_i^* \ (\text{\rm mod }q) \ \text{ for every } \ i = 1, \dots, k^*,
\eeq
then we have by \eqref{eq:6.16}--\eqref{eq:6.17}
\beq
\label{eq:6.19}
\mu_{1,q}^* (p_0) > \mu_{1,q}(p_0)
\eeq
which increases the corresponding term in our crucial sum.
Let $p_0 \notin \mathcal H^*$ be given, and let us fix $i \in [1, k^*]$.
How many different $q$'s do we have at most in $I_1 \subseteq [x/2, x]$ with
\beq
\label{eq:6.20}
p_0 \equiv h_i^* \ (\text{\rm mod } q)?
\eeq
The answer is simple: at most one.
If we had, namely, \eqref{eq:6.20} for $q = q_1, q_2 \in I_1$ $(q_1 \neq q_2)$, then this would imply by $q_1 q_2 \geqslant x^2 / 4 > U \geqslant \max(p_0, h_i^*)$
\beq
\label{eq:6.21}
p_0 \equiv h_i^* \ (\text{\rm mod } q_1 q_2),
\eeq
consequently
\beq
\label{eq:6.22}
p_0 = h_i^* \in \mathcal H^*.
\eeq

So the remaining question is reduced to show that if we delete at most $k^* = O(1)$ terms from the original sum $\sum\limits_q \mu_{1,q}(p_0)$ the sum will be still $\gg \delta \log k$.
But this follows already from the trivial relation $\mu_{m,q}(a) \leqslant 1$, although it is easy to see that even
$\mu_{m,q}(a) \ll_\varepsilon \frac{x^\varepsilon}{q} \ll \frac{1}{x^{1 - \varepsilon}}$
holds for $q \in [x/2, x]$.
($k$ can be chosen sufficiently large compared with $k^*$.)

This completes the proof of Theorems 1', 3', 4' (and thereby Corollary 1').
In case of Theorem 2' we consider a given $N'$ and distinguish the following two cases.
Let us consider a sequence of $99$ exceptional functions.

\smallskip
\emph{Case 1*.}
\[
f_{50} (N') < \log N'(\log_2 N')^{1/2}.
\]

In this case the proof of the original Theorem~\ref{th:2} can be applied to the increasing subset $\left\{f_i(x)\right\}_{i = 1}^{50}$.

\smallskip
\emph{Case 2*.}
\[
f_{50} (N') \geqslant \log N'(\log_2 N')^{1/2}.
\]

In this case the new method of \cite{May2} together with the changes of the present section yields the result for the increasing subset $\left\{f_i(x)\right\}_{i = 50}^{99}$.

Thus, in both cases we obtain a contradiction if we suppose for an increasing sequence of at least $99$ functions that the relation \eqref{eq:1.12} fails.

\begin{remark}
\label{rem:6.1}
With some extra effort it would also be possible to show Theorem 2' with $49$ instead of~$98$.
\end{remark}

\begin{remark}
\label{rem:6.2}
If we define additionally
\beq
\label{eq:6.23}
\lambda_{d_1, \dots, d_k, e_1, \dots, e_k} = 0 \ \text{ if } \ q_0 \mid \prod_{i = 1}^k (d_i e_i),
\eeq
then the whole result can be made effective.
(We remark that actually we have a loss of size $1 + O\left(\frac{1}{q_0}\right) = 1 + o(1)$ due to \eqref{eq:6.23}
but this does not affect the validity of the argument.)
\end{remark}

\bigskip
\noindent
{\bf Acknowledgement:} The author would like to express his sincere
gratitude to Imre Z. Ruzsa, who called his attention that a possible
combination of the methods of Erd\H{o}s--Rankin and Zhang--Maynard--Tao might
lead to stronger results about gaps between consecutive primes.

\noindent
J\'anos Pintz\\
Alfr\'ed R\'enyi Institute of Mathematics,\\
Hungarian Academy of Sciences\\
Budapest, Re\'altanoda u. 13--15\\
H-1053 Hungary\\
e-mail: pintz.janos@renyi.mta.hu

\end{document}